\documentclass[a4paper,12pt]{article}
\usepackage{amssymb}
{}

\newtheorem{thm}{Theorem}[section]

\newcommand{\R}{\mathbb{R}}

\renewcommand{\r}{\rho}

\renewcommand{\a}{\alpha}

\newcommand{\p}{\partial}
\renewcommand{\l}{\lambda}

\newcommand{\bc}{\begin{cor}}
\newcommand{\ec}{\end{cor}}
\newcommand{\bl}{\begin{lem}}
\newcommand{\el}{\end{lem}}
\newcommand{\bp}{\begin{prop}}
\newcommand{\ep}{\end{prop}}
\newcommand{\bt}{\begin{thm}}
\newcommand{\et}{\end{thm}}
\newcommand{\bal}{\begin{array}{ll}}
\newcommand{\ba}{\begin{array}}
\newcommand{\bac}{\begin{array}{ccc}}
\newcommand{\ea}{\end{array}}

\newcommand{\be}{\begin{equation}}
\newcommand{\ee}{\end{equation}}
\newcommand{\pt}{\partial_t}
\newtheorem{lem}[thm]{Lemma}
\newtheorem{cor}[thm]{Corollary}
\newtheorem{prop}[thm]{Proposition}

\begin{document}
\title {\Large\bf{$L^p-L^q$ estimates on the solutions to $u_{tt}-u_{x_1x_1}=\triangle u_t$
  }}

\author { Yongqin Liu{\footnote{email:~yqliu2@yahoo.com.cn}}\\
 {\small\emph{School of Mathematical Sciences, Fudan University, Shanghai, China}}\\
 Yi Zhou{\footnote{email:~yizhou@fudan.ac.cn}}\\
 {\small\emph{School of Mathematical Sciences, Fudan University, Shanghai, China}}\\
 }
\date{}
\maketitle
 \noindent{\bf Abstract}. This paper focuses the study on the $L^p-L^q$ estimates on the solutions
 to an asymmetric wave equation with dissipation which arises in the study for the
 magneto-hydrodynamics by using the method of Green function.
\section{Introduction}

We consider the following Cauchy problem of the asymmetric
dissipative wave equation, \be\label{1a}
 \left\{\begin{array}{ll}
&u_{tt}-u_{x_1x_1}=\triangle u_t, \ \ x=(x_1, x_2, \cdots, x_n)\in \R^n, \ t>0,\\&\\
&u(x, t)|_{t=0}=u_0(x),\\&\\&u_t(x, t)|_{t=0}=u_1(x),
\end{array}
\right. \ee where
$\triangle=\p^2_{x_1}+\p^2_{x_2}+\cdots+\p^2_{x_n},\ n \geq3,$ and
$u_0, u_1$ are given functions.

Equation (\ref{1a}) is an asymmetric wave equation with dissipation
and arises in  the study for the magneto-hydrodynamics (see
\cite{LQ}). Due to the asymmetric structure of (\ref{1a}), the
energy method fails
 in making the $L^p-L^q$ estimates, and new difficulties arise in the
 study for this equation in contrast to that for the symmetric
  semilinear and nonlinear wave equation (see \cite{INZ, LW,LZ,NZ,O1,O2,O3,TY}).
 In this paper, we obtain the $L^p-L^q$ estimates on the solutions to (\ref{1a}) by using
  the method of Green function combined with the technique of
 Fourier analysis.

Notations. We denote generic constants by $C$.
 $L^p(1\leq p\leq\infty)$ is
the usual Lebesgue space with the norm $\|\cdot\|_{L^p}$,
$W^{m,~p},m\in {\mathbb{Z}^+},p\in[1,\infty]$ denotes the usual
Sobolev space with its norm
$$
\|f\|_{W^{m,~p}}:=(\sum\limits_{k=0}\limits^{m}|\partial_x^kf|_p^p)^{1\over
p}.
$$
\section{Green function}
The corresponding Green function to (\ref{1a}) satisfies the
following equation, \be\label{2a}
 \left\{\begin{array}{ll}
&G_{tt}-G_{x_1x_1}=\triangle G_t, \ \ x\in \R^n, \ t>0,\\&\\
&G(x, t)|_{t=0}=0,\\&\\&G_t(x, t)|_{t=0}=\delta(x).
\end{array}
\right. \ee

By Fourier transformation we get that \be\label{2b}
 \left\{\begin{array}{ll}
&\hat{G}_{tt}+|\xi|^2\hat{G}_t+\xi^2_1\hat{G}=0, \ \ \xi=(\xi_1, \xi_2, \cdots, \xi_n)\in \R^n, \ t>0,\\&\\
&\hat{G}(\xi, t)|_{t=0}=0,\\&\\&\hat{G}_t(\xi, t)|_{t=0}=1.
\end{array}
\right. \ee

It yields that $$ \hat{G}(\xi,
t)={{e^{\l_+t}-e^{\l_-t}}\over{\l_0}},
$$
where $\l_{\pm}={{-|\xi|^2\pm\sqrt{|\xi|^4-4\xi_1^2}}\over{2}},
\l_0=\sqrt{|\xi|^4-4\xi_1^2}.$

By Duhamel's principle we know that the solution to (\ref{1a}) can
be expressed as following, \be\label{2c} u(x,t)=(\pt-\triangle)G\ast
u_0(x,t)+G\ast u_1(x,t). \ee

Denote $A=\{\xi\in \R^n; |\xi|^2\leq2|\xi_1|\},\ B=\{\xi\in \R^n;
|\xi|<1\},\ E=\{\xi\in \R^n; |\xi_1|\leq|\xi^{\prime}|\},\
\xi^{\prime}=(\xi_2,\cdots, \xi_n),$ and $ D_r=\{\xi\in \R^n;
|\xi|^2\leq r|\xi_1|\},\ r>2.$ Define $A^c$ the complete set of $A$,
and $ B^c, D_r^c, E^c$ can be defined similarly.

 By direct calculation and induction we have the following
lemma.

\bl\label{21}

 $$ \pt(\pt+|\xi|^2)\hat{G}(\xi,
t)=-\xi_1^2\hat{G}(\xi, t),$$
$$\pt^2(\pt+|\xi|^2)\hat{G}(\xi, t)=-\xi_1^2\hat{G}_t(\xi, t). $$
$\forall h\geq1,$ if $\xi\in B,$ then $$ \pt^{2h}\hat{G}(\xi,
t)=O(|\xi|^{2h})\hat{G}(\xi, t)+O(|\xi|^{2h})\hat{G}_t(\xi, t),
 $$
$$ \pt^{2h+1}\hat{G}(\xi, t)=O(|\xi|^{2h+2})\hat{G}(\xi,
t)+O(|\xi|^{2h})\hat{G}_t(\xi, t), $$
$$ \pt^{2h+1}(\pt+|\xi|^2)\hat{G}(\xi, t)=O(|\xi|^{2h+2})\hat{G}(\xi,
t)+O(|\xi|^{2h+2})\hat{G}_t(\xi, t), $$
$$ \pt^{2h+2}(\pt+|\xi|^2)\hat{G}(\xi, t)=O(|\xi|^{2h+4})\hat{G}(\xi,
t)+O(|\xi|^{2h+2})\hat{G}_t(\xi, t). $$ If $\xi\in B^c,$ then
 $$ \pt^{2h}\hat{G}(\xi,
t)=O(|\xi|^{4h-2})\hat{G}(\xi, t)+O(|\xi|^{4h-2})\hat{G}_t(\xi,
t),$$
$$ \pt^{2h+1}\hat{G}(\xi, t)=O(|\xi|^{4h})\hat{G}(\xi,
t)+O(|\xi|^{4h})\hat{G}_t(\xi, t), $$ $$
\pt^{2h+1}(\pt+|\xi|^2)\hat{G}(\xi, t)=O(|\xi|^{4h})\hat{G}(\xi,
t)+O(|\xi|^{4h})\hat{G}_t(\xi, t), $$
$$ \pt^{2h+2}(\pt+|\xi|^2)\hat{G}(\xi, t)=O(|\xi|^{4h+2})\hat{G}(\xi,
t)+O(|\xi|^{4h+2})\hat{G}_t(\xi, t). $$ \el
 By direct calculation we have the following estimates on $\hat{G}(\xi, t)$ and $\hat{G}_t(\xi, t)$.

 \bl\label{22} If $\xi\in D_{r},\ r>2,$ then
$$ |\hat{G}(\xi, t)|\leq te^{-{{m}\over2}|\xi|^2t},\ |\hat{G}_t(\xi,
t)|\leq (1+|\xi|^2t)e^{-{{m}\over2}|\xi|^2t},$$ here,
$m=1-\sqrt{1-{4\over{{r}^2}}}.$ \el

 {\bf Proof.} If $\xi\in A,$
then $|\xi|^4\leq4\xi_1^2.$ In this case,
$$
\hat{G}(\xi,
t)=e^{-{{|\xi|^2}\over2}t}\cdot{{e^{{{\sqrt{|\xi|^4-4\xi_1^2}}\over2}t}-e^{-{{\sqrt{|\xi|^4-4\xi_1^2}}
\over2}t}}\over{\sqrt{|\xi|^4-4\xi_1^2}}}=te^{-{{|\xi|^2}\over2}t}e^{{{\theta(\xi)\sqrt{|\xi|^4-4\xi_1^2}}\over2}t},
$$ for some $\theta(\xi)\in(-1,1)$. It yields that $|\hat{G}(\xi, t)|\leq te^{-{{|\xi|^2}\over2}t}$.

Similarly, we have that  $$|\hat{G}_t(\xi,
t)|=\left|{{e^{\l_+t}+e^{\l_-t}}\over{2}}-{{|\xi|^2}\over2}\hat{G}(\xi,
t)\right|\leq (1+{{|\xi|^2}}t)e^{-{{|\xi|^2}\over2}t}.$$

If $\xi\in A^c\cap D_{r},$ then $4\xi_1^2<|\xi|^4\leq r^2\xi_1^2.$
In this case, $$ \hat{G}(\xi,
t)=te^{-{{|\xi|^2}\over2}t}e^{{{\theta(\xi)\sqrt{|\xi|^4-4\xi_1^2}}\over2}t}\leq
te^{-{{1-\sqrt{1-{4\over{{r}^2}}}}\over2}|\xi|^2t}.
$$
Similarly, we have that  $$|\hat{G}_t(\xi, t)|\leq
(1+{{|\xi|^2}}t)e^{-{{1-\sqrt{1-{4\over{{r}^2}}}}\over2}|\xi|^2t}.$$
So the lemma is proved. $\Box$

Denote $\hat{G}_1(\xi,t)=\chi(\xi)\hat{G}(\xi,t),$ where $\chi\in
C_0^{\infty}(\R^n),$ ${\rm{supp}}\chi\subset D_{r+1}$ and
$\chi|_{D_r}=1.$  Then we have the following estimates on
$\hat{G}_1(\xi,t).$

\bl\label{23} For any multi-index $\alpha$ and non-negative integer
$l$, we have that,
$$
\bal (1).& \|(\cdot)^{\a}\pt^{l}\hat{G}_1(\cdot,
t)\|_{L^{\infty}}\leq Ct^{-{{|\a|}\over2}-[{{l-1}\over2}]},
\\&\\(2).& \|(\cdot)^{\a}\pt^{l}(\pt+|\cdot|^2)\hat{G}_1(\cdot,
t)\|_{L^{\infty}}\leq Ct^{-{{|\a|}\over2}-[{l\over2}]}, \\&\\(3).&
\|(\cdot)^{\a}\pt^{l}\hat{G}_1(\cdot, t)\|_{L^{2}}\leq
Ct^{-{{|\a|}\over2}-[{{l-1}\over2}]-{n\over4}}, \\&\\(4).&
\|(\cdot)^{\a}\pt^{l}(\pt+|\cdot|^2)\hat{G}_1(\cdot,
t)\|_{L^{2}}\leq Ct^{-{{|\a|}\over2}-[{l\over2}]-{n\over4}},
\\&\\(5).& \|(\cdot)^{\a}\pt^{l}\hat{G}_1(\cdot, t)\|_{L^{1}}\leq
Ct^{-{{|\a|}\over2}-[{{l-1}\over2}]-{n\over2}}, \\&\\(6).&
\|(\cdot)^{\a}\pt^{l}(\pt+|\cdot|^2)\hat{G}_1(\cdot,
t)\|_{L^{1}}\leq Ct^{-{{|\a|}\over2}-[{l\over2}]-{n\over2}}, \ea
$$ where $[\cdot]$ is Gauss' symbol.
\el

{\bf Proof.} By using lemma \ref{21} and lemma \ref{22} we give the
proof.

(1). As $l=0,$ $$|\xi^{\alpha}\chi(\xi)\hat{G}(\xi,t)|\leq
C|\xi|^{|\alpha|}te^{-{{m}\over2}|\xi|^2t}\leq
Ct^{-{{|\alpha|}\over2}+1}.$$

As $l=1,$ $$|\xi^{\alpha}\chi(\xi)\pt\hat{G}(\xi,t)|\leq
C|\xi|^{|\alpha|}(1+|\xi|^2t)e^{-{{m}\over2}|\xi|^2t}\leq
Ct^{-{{|\alpha|}\over2}}.$$
 For $h \geq1,$ as $l=2h,$
$$\bal|\xi^{\alpha}\chi(\xi)\pt^{2h}\hat{G}(\xi,t)|_B&=|\xi^{\alpha}\chi(\xi)\{O(|\xi|^{2h})\hat{G}(\xi,t)
+O(|\xi|^{2h})\hat{G}_t(\xi,t)\}|\\&\\&\leq
C|\xi|^{|\alpha|+2h}(t+1+|\xi|^2t)e^{-{{m}\over2}|\xi|^2t}\leq
Ct^{-{{|\alpha|}\over2}-h+1},
 \ea$$
 $$\bal|\xi^{\alpha}\chi(\xi)\pt^{2h}\hat{G}(\xi,t)|_{B^c}&=|\xi^{\alpha}\chi(\xi)\{O(|\xi|^{4h-2})\hat{G}(\xi,t)
+O(|\xi|^{4h-2})\hat{G}_t(\xi,t)\}|\\&\\&\leq
Ct^{-{{|\alpha|}\over2}-2h+2}\leq Ct^{-{{|\alpha|}\over2}-h+1};\ea$$
 as $l=2h+1,$
$$\bal|\xi^{\alpha}\chi(\xi)\pt^{2h+1}\hat{G}(\xi,t)|_B&=|\xi^{\alpha}\chi(\xi)\{O(|\xi|^{2h+2})\hat{G}(\xi,t)
+O(|\xi|^{2h})\hat{G}_t(\xi,t)\}|\\&\\&\leq
Ct^{-{{|\alpha|}\over2}-h},\ea$$
$$\bal|\xi^{\alpha}\chi(\xi)\pt^{2h+1}\hat{G}(\xi,t)|_{B^c}&=|\xi^{\alpha}\chi(\xi)\{O(|\xi|^{4h})\hat{G}(\xi,t)
+O(|\xi|^{4h})\hat{G}_t(\xi,t)\}|\\&\\&\leq
Ct^{-{{|\alpha|}\over2}-2h+1}\leq Ct^{-{{|\alpha|}\over2}-h}.\ea$$
Thus (1) is proved.

(2).  As $l=0,$
$$|\xi^{\alpha}\chi(\xi)(\pt+|\xi|^2)\hat{G}(\xi,t)|\leq
C|\xi|^{|\alpha|}(1+2|\xi|^2t)e^{-{{m}\over2}|\xi|^2t}\leq
Ct^{-{{|\alpha|}\over2}}.$$ As $l=1,$
$$|\xi^{\alpha}\chi(\xi)\pt(\pt+|\xi|^2)\hat{G}(\xi,t)|=|\xi^{\a}\chi(\xi)\xi_1^2\hat{G}(\xi,t)|\leq
Ct^{-{{|\alpha|}\over2}}.$$ As $l=2,$
$$|\xi^{\alpha}\chi(\xi)\pt^2(\pt+|\xi|^2)\hat{G}(\xi,t)|=|\xi^{\a}\chi(\xi)\xi_1^2\hat{G}_t(\xi,t)|\leq
Ct^{-{{|\alpha|}\over2}-1}.$$
 For $h \geq1,$ as $l=2h+1,$
$$\bal&\ |\xi^{\alpha}\chi(\xi)\pt^{2h+1}(\pt+|\xi|^2)\hat{G}(\xi,t)|_B\\&\\&=|\xi^{\alpha}\chi(\xi)\{O(|\xi|^{2h+2})\hat{G}(\xi,t)
+O(|\xi|^{2h+2})\hat{G}_t(\xi,t)\}|\\&\\&\leq
Ct^{-{{|\alpha|}\over2}-h},
 \ea$$
 $$\bal&\ |\xi^{\alpha}\chi(\xi)\pt^{2h+1}(\pt+|\xi|^2)\hat{G}(\xi,t)|_{B^c}\\&\\&=|\xi^{\alpha}\chi(\xi)\{O(|\xi|^{4h})\hat{G}(\xi,t)
+O(|\xi|^{4h})\hat{G}_t(\xi,t)\}|\\&\\&\leq
Ct^{-{{|\alpha|}\over2}-2h+1}\leq Ct^{-{{|\alpha|}\over2}-h};\ea$$
as $l=2h+2,$
$$\bal&\ |\xi^{\alpha}\chi(\xi)\pt^{2h+2}(\pt+|\xi|^2)\hat{G}(\xi,t)|_B\\&\\&=|\xi^{\alpha}\chi(\xi)\{O(|\xi|^{2h+4})\hat{G}(\xi,t)
+O(|\xi|^{2h+2})\hat{G}_t(\xi,t)\}|\\&\\&\leq
Ct^{-{{|\alpha|}\over2}-h-1},
 \ea$$
 $$\bal&\ |\xi^{\alpha}\chi(\xi)\pt^{2h+2}(\pt+|\xi|^2)\hat{G}(\xi,t)|_{B^c}\\&\\&=|\xi^{\alpha}\chi(\xi)\{O(|\xi|^{4h+2})\hat{G}(\xi,t)
+O(|\xi|^{4h+2})\hat{G}_t(\xi,t)\}|\\&\\&\leq
Ct^{-{{|\alpha|}\over2}-2h}\leq Ct^{-{{|\alpha|}\over2}-h-1}.\ea$$
Thus (2) is proved.

(3). As $l=0,$ $$ \bal \|(\cdot)^{\a}\hat{G}_1(\cdot,
t)\|_{L^2}&=(\int_{\R^n}|\xi^{\a}\chi(\xi)\hat{G}(\xi,
t)|^2d\xi)^{1\over2}\\&\\&\leq
C(\int_{\R^n}|\xi|^{2|\a|}t^2e^{-m|\xi|^2t}d\xi)^{1\over2}\\&\\&\leq
Ct^{-{{|\alpha|}\over2}+1-{n\over4}}. \ea
$$
As $l=1,$ $$ \bal \|(\cdot)^{\a}\pt\hat{G}_1(\cdot,
t)\|_{L^2}&=(\int_{\R^n}|\xi^{\a}\chi(\xi)\pt\hat{G}(\xi,
t)|^2d\xi)^{1\over2}\\&\\&\leq
C(\int_{\R^n}|\xi|^{2|\a|}(1+|\xi|^2t)^2e^{-m|\xi|^2t}d\xi)^{1\over2}\\&\\&\leq
Ct^{-{{|\alpha|}\over2}-{n\over4}}. \ea
$$
For $h\geq1,$ as $l=2h,$ $$ \bal&\
\|(\cdot)^{\a}\pt^{2h}\hat{G}_1(\cdot,
t)\|_{L^2}\\&\\&=(\int_{B}|\xi^{\a}\chi(\xi)\{O(|\xi|^{2h})\hat{G}(\xi,t)
+O(|\xi|^{2h})\hat{G}_t(\xi,t)\}|^2d\xi\\&\\&+\int_{B^c}|\xi^{\a}\chi(\xi)\{O(|\xi|^{4h-2})\hat{G}(\xi,t)
+O(|\xi|^{4h-2})\hat{G}_t(\xi,t)\}|^2d\xi)^{1\over2}\\&\\&\leq
Ct^{-{{|\alpha|}\over2}-h+1-{n\over4}}; \ea
$$
as $l=2h+1,$ $$ \bal &\ \|(\cdot)^{\a}\pt^{2h+1}\hat{G}_1(\cdot,
t)\|_{L^2}\\&\\&=(\int_{B}|\xi^{\a}\chi(\xi)\{O(|\xi|^{2h+2})\hat{G}(\xi,t)
+O(|\xi|^{2h})\hat{G}_t(\xi,t)\}|^2d\xi\\&\\&+\int_{B^c}|\xi^{\a}\chi(\xi)\{O(|\xi|^{4h})\hat{G}(\xi,t)
+O(|\xi|^{4h})\hat{G}_t(\xi,t)\}|^2d\xi)^{1\over2}\\&\\&\leq
Ct^{-{{|\alpha|}\over2}-h-{n\over4}}. \ea
$$
Thus (3) is proved.

(4). As $l=0,$ $$ \bal \|(\cdot)^{\a}(\pt+|\cdot|^2)\hat{G}_1(\cdot,
t)\|_{L^2}&=(\int_{\R^n}|\xi^{\a}\chi(\xi)(\pt+|\xi|^2)\hat{G}(\xi,
t)|^2d\xi)^{1\over2}\\&\\&\leq Ct^{-{{|\alpha|}\over2}-{n\over4}}.
\ea
$$
As $l=1,$ $$ \bal \|(\cdot)^{\a}\pt(\pt+|\cdot|^2)\hat{G}_1(\cdot,
t)\|_{L^2}&=(\int_{\R^n}|\xi^{\a}\chi(\xi)\pt(\pt+|\xi|^2)\hat{G}(\xi,
t)|^2d\xi)^{1\over2}\\&\\&=(\int_{\R^n}|\xi^{\a}\chi(\xi)\xi_1^2\hat{G}(\xi,
t)|^2d\xi)^{1\over2}\\&\\&\leq Ct^{-{{|\alpha|}\over2}-{n\over4}}.
\ea
$$
As $l=2,$ $$ \bal \|(\cdot)^{\a}\pt^2(\pt+|\cdot|^2)\hat{G}_1(\cdot,
t)\|_{L^2}&=(\int_{\R^n}|\xi^{\a}\chi(\xi)\pt^2(\pt+|\xi|^2)\hat{G}(\xi,
t)|^2d\xi)^{1\over2}\\&\\&=(\int_{\R^n}|\xi^{\a}\chi(\xi)\xi_1^2\hat{G}_t(\xi,
t)|^2d\xi)^{1\over2}\\&\\&\leq Ct^{-{{|\alpha|}\over2}-1-{n\over4}}.
\ea
$$
For $h\geq1,$ as $l=2h+1,$ $$ \bal &\
\|(\cdot)^{\a}\pt^{2h+1}(\pt+|\cdot|^2)\hat{G}_1(\cdot,
t)\|_{L^2}\\&\\&=(\int_{B}|\xi^{\a}\chi(\xi)\{O(|\xi|^{2h+2})\hat{G}(\xi,t)
+O(|\xi|^{2h+2})\hat{G}_t(\xi,t)\}|^2d\xi\\&\\&+\int_{B^c}|\xi^{\a}\chi(\xi)\{O(|\xi|^{4h})\hat{G}(\xi,t)
+O(|\xi|^{4h})\hat{G}_t(\xi,t)\}|^2d\xi)^{1\over2}\\&\\&\leq
Ct^{-{{|\alpha|}\over2}-h-{n\over4}}; \ea
$$
as $l=2h+2,$ $$ \bal &\
\|(\cdot)^{\a}\pt^{2h+2}(\pt+|\cdot|^2)\hat{G}_1(\cdot,
t)\|_{L^2}\\&\\&=(\int_{B}|\xi^{\a}\chi(\xi)\{O(|\xi|^{2h+4})\hat{G}(\xi,t)
+O(|\xi|^{2h+2})\hat{G}_t(\xi,t)\}|^2d\xi\\&\\&+\int_{B^c}|\xi^{\a}\chi(\xi)\{O(|\xi|^{4h+2})\hat{G}(\xi,t)
+O(|\xi|^{4h+2})\hat{G}_t(\xi,t)\}|^2d\xi)^{1\over2}\\&\\&\leq
Ct^{-{{|\alpha|}\over2}-h-1-{n\over4}}. \ea
$$
Thus (4) is proved.

(5). As $l=0,$ $$ \bal \|(\cdot)^{\a}\hat{G}_1(\cdot,
t)\|_{L^1}&=\int_{\R^n}|\xi^{\a}\chi(\xi)\hat{G}(\xi,
t)|d\xi\\&\\&\leq Ct^{-{{|\alpha|}\over2}+1-{n\over2}}. \ea
$$
As $l=1,$ $$ \bal \|(\cdot)^{\a}\pt\hat{G}_1(\cdot,
t)\|_{L^1}&=\int_{\R^n}|\xi^{\a}\chi(\xi)\pt\hat{G}(\xi,
t)|d\xi\\&\\&\leq Ct^{-{{|\alpha|}\over2}-{n\over2}}. \ea
$$
For $h\geq1,$ as $l=2h,$ $$ \bal &\
\|(\cdot)^{\a}\pt^{2h}\hat{G}_1(\cdot,
t)\|_{L^1}\\&\\&=\int_{B}|\xi^{\a}\chi(\xi)\{O(|\xi|^{2h})\hat{G}(\xi,t)
+O(|\xi|^{2h})\hat{G}_t(\xi,t)\}|d\xi\\&\\&+\int_{B^c}|\xi^{\a}\chi(\xi)\{O(|\xi|^{4h-2})\hat{G}(\xi,t)
+O(|\xi|^{4h-2})\hat{G}_t(\xi,t)\}|d\xi\\&\\&\leq
Ct^{-{{|\alpha|}\over2}-h+1-{n\over2}}; \ea
$$
as $l=2h+1,$ $$ \bal &\ \|(\cdot)^{\a}\pt^{2h+1}\hat{G}_1(\cdot,
t)\|_{L^1}\\&\\&=\int_{B}|\xi^{\a}\chi(\xi)\{O(|\xi|^{2h+2})\hat{G}(\xi,t)
+O(|\xi|^{2h})\hat{G}_t(\xi,t)\}|d\xi\\&\\&+\int_{B^c}|\xi^{\a}\chi(\xi)\{O(|\xi|^{4h})\hat{G}(\xi,t)
+O(|\xi|^{4h})\hat{G}_t(\xi,t)\}|d\xi\\&\\&\leq
Ct^{-{{|\alpha|}\over2}-h-{n\over2}}. \ea
$$
Thus (5) is proved.

(6). As $l=0,$ $$ \bal \|(\cdot)^{\a}(\pt+|\cdot|^2)\hat{G}_1(\cdot,
t)\|_{L^1}&=\int_{\R^n}|\xi^{\a}\chi(\xi)(\pt+|\xi|^2)\hat{G}(\xi,
t)|d\xi\\&\\&\leq Ct^{-{{|\alpha|}\over2}-{n\over2}}. \ea
$$
As $l=1,$ $$ \bal \|(\cdot)^{\a}\pt(\pt+|\cdot|^2)\hat{G}_1(\cdot,
t)\|_{L^1}&=\int_{\R^n}|\xi^{\a}\chi(\xi)\pt(\pt+|\xi|^2)\hat{G}(\xi,
t)|^2d\xi\\&\\&=\int_{\R^n}|\xi^{\a}\chi(\xi)\xi_1^2\hat{G}(\xi,
t)|d\xi\\&\\&\leq Ct^{-{{|\alpha|}\over2}-{n\over2}}. \ea
$$
As $l=2,$ $$ \bal \|(\cdot)^{\a}\pt^2(\pt+|\cdot|^2)\hat{G}_1(\cdot,
t)\|_{L^1}&=\int_{\R^n}|\xi^{\a}\chi(\xi)\pt^2(\pt+|\xi|^2)\hat{G}(\xi,
t)|d\xi\\&\\&=\int_{\R^n}|\xi^{\a}\chi(\xi)\xi_1^2\hat{G}_t(\xi,
t)|d\xi\\&\\&\leq Ct^{-{{|\alpha|}\over2}-1-{n\over2}}. \ea
$$
For $h\geq1,$ as $l=2h+1,$ $$ \bal &\
\|(\cdot)^{\a}\pt^{2h+1}(\pt+|\cdot|^2)\hat{G}_1(\cdot,
t)\|_{L^1}\\&\\&=\int_{B}|\xi^{\a}\chi(\xi)\{O(|\xi|^{2h+2})\hat{G}(\xi,t)
+O(|\xi|^{2h+2})\hat{G}_t(\xi,t)\}|d\xi\\&\\&+\int_{B^c}|\xi^{\a}\chi(\xi)\{O(|\xi|^{4h})\hat{G}(\xi,t)
+O(|\xi|^{4h})\hat{G}_t(\xi,t)\}|d\xi)\\&\\&\leq
Ct^{-{{|\alpha|}\over2}-h-{n\over2}}; \ea
$$
as $l=2h+2,$ $$ \bal &\
\|(\cdot)^{\a}\pt^{2h+2}(\pt+|\cdot|^2)\hat{G}_1(\cdot,
t)\|_{L^1}\\&\\&=\int_{B}|\xi^{\a}\chi(\xi)\{O(|\xi|^{2h+4})\hat{G}(\xi,t)
+O(|\xi|^{2h+2})\hat{G}_t(\xi,t)\}|d\xi\\&\\&+\int_{B^c}|\xi^{\a}\chi(\xi)\{O(|\xi|^{4h+2})\hat{G}(\xi,t)
+O(|\xi|^{4h+2})\hat{G}_t(\xi,t)\}|d\xi\\&\\&\leq
Ct^{-{{|\alpha|}\over2}-h-1-{n\over2}}. \ea
$$
Thus (6) is proved. $\Box$
\section{$L^p-L^q$ estimates}
 The first result is the following theorem.

 \bt\label{31}
Assume that $u_0,\ u_1\in L^p(\R^n), n\geq3,$ and there exists a
constant $r>2,$ such that ${\rm{supp}}\hat{u_0}, \
{\rm{supp}}\hat{u_1}\subset D_r,$ then for any multi-index $\a$ and
non-negative integer $l$, we have that,

(1). If $p\in [1,2)$, it holds that, $\forall q\in
[{{2p}\over{2-p}}, +\infty],$
$$\|\partial_x^{\a}\pt^lu(\cdot,t)\|_{L^q}\leq Ct^{-{{|\a|}\over2}-[{{l}\over2}]-{{n}\over2}({1\over p}-{1\over
q})}\|u_0\|_{L^p}+Ct^{-{{|\a|}\over2}-[{{l-1}\over2}]-{{n}\over2}({1\over
p}-{1\over q})}\|u_1\|_{L^p};$$

(2). If $p=2,$ it holds that, $\forall q\in [2, +\infty],$
$$
\|\partial_x^{\a}\pt^lu(\cdot,t)\|_{L^q}\leq
Ct^{-{{|\a|}\over2}-[{{l}\over2}]-{{n}\over2}({1\over 2}-{1\over
q})}\|u_0\|_{L^2}+Ct^{-{{|\a|}\over2}-[{{l-1}\over2}]-{{n}\over2}({1\over
2}-{1\over q})}\|u_1\|_{L^2}.
$$ where $\|\cdot\|_{L^p}=\|\cdot\|_{L^p(\R^n)}.$
 \et
{\bf Proof.}
 Since
${\rm{supp}}\hat{u_0},\ {\rm{supp}}\hat{u_1}\subset D_r,$ from
(\ref{2c}) we have that
$$\bal\hat{u}(\xi,t)&=(\pt+|\xi|^2)\hat{G}(\xi,t)\hat{u}_0(\xi)+\hat{G}(\xi,t)\hat{u}_1(\xi)\\&\\&=
(\pt+|\xi|^2)\hat{G}_1(\xi,t)\hat{u}_0(\xi)+\hat{G}_1(\xi,t)\hat{u}_1(\xi),\ea$$
it yields that
$$
u(x,t)=(\pt-\triangle)G_1\ast u_0(x,t)+G_1\ast u_1(x,t).
$$
(1). Denote $1+{1\over q}={1\over p}+{1\over \r}.$ If
$q\geq{{2p}\over{2-p}},$ then $\r\geq2.$ In view of lemma
{\ref{23}}, Young inequality and interpolation formula, we have
that,
$$
\bal&\|\p_x^{\a}\pt^l u(\cdot, t)\|_{L^q}\\&\\&\leq \|\p_x^{\a}\pt^l
(\pt-\triangle)G_1\ast u_0(\cdot, t)\|_{L^q}+\|\p_x^{\a}\pt^l
G_1\ast u_1(\cdot, t)\|_{L^q}\\&\\&\leq\|\p_x^{\a}\pt^l
(\pt-\triangle)G_1(\cdot,
t)\|_{L^{\r}}\|u_0\|_{L^p}+\|\p_x^{\a}\pt^l G_1(\cdot,
t)\|_{L^{\r}}\|u_1\|_{L^p}\\&\\&\leq \|\p_x^{\a}\pt^l
(\pt-\triangle)G_1(\cdot, t)\|_{L^{2}}^{2\over{\r}}\|\p_x^{\a}\pt^l
(\pt-\triangle)G_1(\cdot,
t)\|_{L^{\infty}}^{1-{2\over{\r}}}\|u_0\|_{L^p}\\&\\&\ \
+\|\p_x^{\a}\pt^l G_1(\cdot,
t)\|_{L^{2}}^{2\over{\r}}\|\p_x^{\a}\pt^l G_1(\cdot,
t)\|_{L^{\infty}}^{1-{2\over{\r}}}\|u_1\|_{L^p}\\&\\&\leq C
\|(\cdot)^{\a}\pt^l(\pt+|\cdot|^2)\hat{G}_1(\cdot,
t)\|_{L^2}^{2\over{\r}}\|(\cdot)^{\a}\pt^l(\pt+|\cdot|^2)\hat{G}_1(\cdot,
t)\|_{L^1}^{1-{2\over {\r}}}\|u_0\|_{L^p}\\&\\&\ \ +C
\|(\cdot)^{\a}\pt^l\hat{G}_1(\cdot,
t)\|_{L^2}^{2\over{\r}}\|(\cdot)^{\a}\pt^l\hat{G}_1(\cdot,
t)\|_{L^1}^{1-{2\over {\r}}}\|u_1\|_{L^p}\\&\\&\leq
Ct^{-{{|\a|}\over2}-[{{l}\over2}]-{{n}\over2}(1-{1\over{\r}})}\|u_0\|_{L^p}+
Ct^{-{{|\a|}\over2}-[{{l-1}\over2}]-{{n}\over2}(1-{1\over{\r}})}\|u_1\|_{L^p}\\&\\&=
Ct^{-{{|\a|}\over2}-[{{l}\over2}]-{{n}\over2}({1\over p}-{1\over
q})}\|u_0\|_{L^p}+Ct^{-{{|\a|}\over2}-[{{l-1}\over2}]-{{n}\over2}({1\over
p}-{1\over q})}\|u_1\|_{L^p}.
 \ea
$$
Thus (1) is proved.

(2).  By using lemma {\ref{23}} and Young inequality, we have that
$$
\bal&\|\p_x^{\a}\pt^l u(\cdot, t)\|_{L^{\infty}}\\&\\&\leq
\|\p_x^{\a}\pt^l (\pt-\triangle)G_1\ast u_0(\cdot,
t)\|_{L^{\infty}}+\|\p_x^{\a}\pt^l G_1\ast u_1(\cdot,
t)\|_{L^{\infty}}\\&\\&\leq\|\p_x^{\a}\pt^l
(\pt-\triangle)G_1(\cdot, t)\|_{L^{2}}\|u_0\|_{L^2}+\|\p_x^{\a}\pt^l
G_1(\cdot, t)\|_{L^{2}}\|u_1\|_{L^2}\\&\\&=
\|(\cdot)^{\a}\pt^l(\pt+|\cdot|^2)\hat{G}_1(\cdot,
t)\|_{L^2}\|u_0\|_{L^2}+\|(\cdot)^{\a}\pt^l\hat{G}_1(\cdot,
t)\|_{L^2}\|u_1\|_{L^2}\\&\\&\leq
Ct^{-{{|\a|}\over2}-[{{l}\over2}]-{{n}\over4}}\|u_0\|_{L^2}+
Ct^{-{{|\a|}\over2}-[{{l-1}\over2}]-{{n}\over4}}\|u_1\|_{L^2}.
 \ea
$$
$$
\bal&\|\p_x^{\a}\pt^l u(\cdot, t)\|_{L^{2}}\\&\\&\leq
\|\p_x^{\a}\pt^l (\pt-\triangle)G_1\ast u_0(\cdot,
t)\|_{L^2}+\|\p_x^{\a}\pt^l G_1\ast u_1(\cdot, t)\|_{L^2}\\&\\&=
\|(\cdot)^{\a}\pt^l(\pt+|\cdot|^2)\hat{G}_1(\cdot,
t)\hat{u_0}\|_{L^2}+\|(\cdot)^{\a}\pt^l\hat{G}_1(\cdot,
t)\hat{u_1}\|_{L^2}\\&\\&\leq
\|(\cdot)^{\a}\pt^l(\pt+|\cdot|^2)\hat{G}_1(\cdot,
t)\|_{L^{\infty}}\|u_0\|_{L^2}+\|(\cdot)^{\a}\pt^l\hat{G}_1(\cdot,
t)\|_{L^{\infty}}\|u_1\|_{L^2}\\&\\&\leq
Ct^{-{{|\a|}\over2}-[{{l}\over2}]}\|u_0\|_{L^2}+
Ct^{-{{|\a|}\over2}-[{{l-1}\over2}]}\|u_1\|_{L^2}.
 \ea
$$

By using interpolation formula, we have that, for $q\in[2, \infty],$
$$
\bal \|\p_x^{\a}\pt^l u(\cdot, t)\|_{L^q}&\leq \|\p_x^{\a}\pt^l
u(\cdot, t)\|_{L^2}^{2\over q}\|\p_x^{\a}\pt^l u(\cdot,
t)\|_{L^{\infty}}^{1-{2\over q}}\\&\\&\leq
Ct^{-{{|\a|}\over2}-[{{l}\over2}]-{{n}\over2}({1\over 2}-{1\over
q})}\|u_0\|_{L^2}+Ct^{-{{|\a|}\over2}-[{{l-1}\over2}]-{{n}\over2}({1\over
2}-{1\over q})}\|u_1\|_{L^2}. \ea
$$
 Thus (2) is proved and so the theorem is proved.~~~$\Box$\\
\section{$L^1-L^{\infty}$ estimates}
In the last section we have some good results by assuming the
Fourier transform of the initial data have special compact support.
In this section we intend to obtain the $L^1-L^{\infty}$ decay
estimates without the assumption of compact support.

Denote $\hat{G}_2(\xi,t)=(1-\chi(\xi))\hat{G}(\xi,t),$ then
$\hat{G}(\xi,t)=\hat{G}_1(\xi,t)+\hat{G}_2(\xi,t)$. Denote
$u(x,t)=v(x,t)+w(x,t),$ where $v(x,t)$ and $w(x,t)$ satisfying
$$\hat{v}(\xi,t)=(\pt+|\xi|^2)\hat{G}_1(\xi,t)\hat{u}_0(\xi)+\hat{G}_1(\xi,t)\hat{u}_1(\xi)$$
$$
\hat{w}(\xi,t)=(\pt+|\xi|^2)\hat{G}_2(\xi,t)\hat{u}_0(\xi)+\hat{G}_2(\xi,t)\hat{u}_1(\xi)
$$

By direct calculation we have that
$$
\hat{w}=(1-\chi(\xi))[{1\over2}(e^{\lambda_+t}+e^{\lambda_-t}+|\xi|^2\hat{G})\hat{u}_0+\hat{G}\hat{u}_1].
$$

 Since $\chi\in
C_0^{\infty}(\R^n),$ ${\rm{supp}}\chi\subset D_{r+1}$ and
$\chi|_{D_r}=1,$ we have that
$$
\begin{array}{ll}
\|\hat{w}(\cdot, t)\|_{L^1(\R^n)}&\leq
\int_{D_r^c}[{1\over2}(e^{\lambda_+t}+e^{\lambda_-t}+|\xi|^2\hat{G}(\xi,t))|\hat{u}_0(\xi)|+\hat{G}(\xi,t)|\hat{u}_1(\xi)
|]d\xi
\\&\\&\leq{1\over2}\int_{D_r^c\cap E}e^{\lambda_+t}|\hat{u}_0(\xi)|d\xi+{1\over2}\int_{D_r^c\cap E^c}
e^{\lambda_+t}|\hat{u}_0(\xi)|d\xi\\&\\&\ \
+{1\over2}\int_{D_r^c}e^{\lambda_-t}|\hat{u}_0(\xi)|d\xi+{1\over2}\int_{D_r^c}|\xi|^2\hat{G}(\xi,t)|\hat{u}_0(\xi)|d\xi
\\&\\&\ \ +\int_{D_r^c\cap
E}\hat{G}(\xi,t)|\hat{u}_1(\xi)|d\xi+\int_{D_r^c\cap
E^c}\hat{G}(\xi,t)|\hat{u}_1(\xi)|d\xi.\end{array}
$$

Next we come to estimate the six terms respectively, where we will
use the facts that $\lambda_+\leq -{{\xi_1^2}\over{|\xi|^2}}$,
 if $\xi\in D_r^c.$

 If $\xi\in E$, then $|\xi|^2\leq2|\xi^{\prime}|^2.$
$$\begin{array}{ll}&
\int_{D_r^c\cap
E}e^{\lambda_+t}|\hat{u}_0(\xi)|d\xi\leq\int_{D_r^c\cap
E}e^{-{{\xi_1^2}\over{|\xi^{\prime}|^2}}t}|\hat{u}_0(\xi)|d\xi\\&\\&\leq\int_{D_r^c\cap
E}e^{-{{\xi_1^2}\over{|\xi^{\prime}|^2}}t}{1\over{|\xi^{\prime}|^{n-1}(1+|\xi^{\prime}|^2)}}d\xi\sup\limits_{\xi\in
\R^n}\{|\xi^{\prime}|^{n-1}(1+|\xi^{\prime}|^2)|\hat{u}_0(\xi)|\}\\&\\&\leq
Ct^{-{1\over2}}\|u_0\|_{W^{n+1,1}(\R^n)}.
 \end{array}$$
 If $\xi\in E^c$, then $|\xi|^2\leq2|\xi_1|^2.$
 $$\begin{array}{ll}&
\int_{D_r^c\cap
E^c}e^{\lambda_+t}|\hat{u}_0(\xi)|d\xi\leq\int_{D_r^c\cap
E^c}e^{-{t\over2}}|\hat{u}_0(\xi)|d\xi\\&\\&\leq
e^{-{{t}\over2}}\int_{D_r^c\cap
E^c}{1\over{(1+|\xi|^2)^{{n+1}\over2}}}d\xi\sup\limits_{\xi\in
\R^n}\{(1+|\xi|^2)^{{{n+1}\over2}}|\hat{u}_0(\xi)|\}\\&\\&\leq
Ce^{-{{t}\over2}}\|u_0\|_{W^{n+1,1}(\R^n)}.
 \end{array}$$
 Since $\lambda_-\leq e^{-{{|\xi|^2}\over2}t}, \forall \xi\in D_r^c,$
$$
\int_{D_r^c}e^{\lambda_-t}|\hat{u}_0(\xi)|d\xi\leq\int_{D_r^c}e^{-{{|\xi|^2}\over2}t}|\hat{u}_0(\xi)|d\xi
\leq Ct^{-{n\over2}}\|u_0\|_{L^1(\R^n)}.
$$
Since $|\xi|^2\hat{G}\leq Ce^{\lambda_+t},\forall \xi\in D_r^c,$ by
the previous calculation we have that
$$
\int_{D_r^c}|\xi|^2\hat{G}(\xi,t)|\hat{u}_0(\xi)|d\xi\leq
C\int_{D_r^c}e^{\lambda_+t}|\hat{u}_0(\xi)|d\xi\leq
Ct^{-{1\over2}}\|u_0\|_{W^{n+1,1}(\R^n)}.
$$
Since $\hat{G}\leq{{2\lambda_+}\over{\lambda_0}}, \forall \xi\in
D_r^c,$
$$\begin{array}{ll}
&\int_{D_r^c\cap E}\hat{G}(\xi,t)|\hat{u}_1(\xi)|d\xi\\&\\&\leq
\int_{D_r^c\cap
E}{{2e^{\lambda_+t}}\over{\lambda_0}}|\hat{u}_1(\xi)|d\xi\leq
C\int_{D_r^c\cap
E}e^{-{{\xi_1^2}\over{|\xi^{\prime}|^2}}t}{1\over{|\xi^{\prime}|^2}}|\hat{u}_1(\xi)|d\xi\\&\\&\leq
C\int_{D_r^c\cap
E}e^{-{{\xi_1^2}\over{|\xi^{\prime}|^2}}t}{1\over{|\xi^{\prime}|^{n-1}(1+|\xi^{\prime}|^2)}}d\xi
\sup\limits_{\xi\in
\R^n}\{|\xi^{\prime}|^{n-3}(1+|\xi^{\prime}|^2)|\hat{u}_1(\xi)|\}
\\&\\&\leq Ct^{-{1\over2}}\|u_0\|_{W^{n-1,1}(\R^n)}.
\end{array}$$
Since $\hat{G}\leq te^{\lambda_+t}, \forall \xi\in D_r^c$,
$$\begin{array}{ll}
&\int_{D_r^c\cap E^c}\hat{G}(\xi,t)|\hat{u}_1(\xi)|d\xi\leq
\int_{D_r^c\cap E}te^{\lambda_+t}|\hat{u}_1(\xi)|d\xi\\&\\&\leq
\int_{D_r^c\cap E}te^{-{t\over2}}|\hat{u}_1(\xi)|d\xi\leq
Ce^{-{{t}\over3}}\|u_0\|_{W^{n+1,1}(\R^n)}.
\end{array}$$

By combining the six inequalities, we obtain that
$$
\|\hat{w}(\cdot,t)\|_{L^1(\R^n)}\leq
Ct^{-{1\over2}}\|(u_0,u_1)\|_{W^{n+1,1}(\R^n)}.
$$
It yields that
$$
\|w(\cdot,t)\|_{L^{\infty}(\R^n)}\leq
Ct^{-{1\over2}}\|(u_0,u_1)\|_{W^{n+1,1}(\R^n)}.
$$

From theorem \ref{31} we know that $v(x,t)$ satisfies,
$$
\|v(\cdot,t)\|_{L^{\infty}(\R^n)}\leq
Ct^{-({n\over2}-1)}\|(u_0,u_1)\|_{L^1(\R^n)}.
$$
Since $n\geq3,$ the solution $u(x,t)$ to the equation (\ref{1a})
satisfies,
$$
\|u(\cdot,t)\|_{L^{\infty}(\R^n)}\leq
Ct^{-{1\over2}}\|(u_0,u_1)\|_{W^{n+1,1}(\R^n)}.
$$

Thus we obtain the following decay estimate. \bt\label{41} Assume
that $u_0,\ u_1\in W^{n+1,1}(\R^n), n\geq3,$ then the solution
$u(x,t)$ to the equation (\ref{1a}) satisfies,
$$
\|u(\cdot,t)\|_{L^{\infty}(\R^n)}\leq
Ct^{-{1\over2}}\|(u_0,u_1)\|_{W^{n+1,1}(\R^n)}.
$$
\et


\end{document}